\documentclass{article}

\usepackage{amssymb, amsmath, amsthm}

\theoremstyle{plain}
\newtheorem{thm}{Theorem}[section]
\newtheorem{lemma}[thm]{Lemma}
\newtheorem{corollary}[thm]{Corollary}
\newtheorem{proposition}[thm]{Proposition}

\theoremstyle{definition}

\newtheorem{remark}[thm]{Remark}

\newcommand{\N}{\mathbb N}
\newcommand{\DD}{\mathbb D}
\newcommand{\C}{\mathbb C}

\newcommand{\Z}{\mathbb Z}

\newcommand{\T}{\mathbb T}

\newcommand{\om}{\Omega}

\usepackage{mathrsfs}

\newcommand{\ind}{\textnormal{ind}}

\title{\bf \center{Generic boundary behaviour of Taylor series in Hardy and Bergman spaces}}
\author{Hans-Peter Beise and J\"urgen M\"uller\footnote{University of Trier, FB IV, Mathematics, D-54286 Trier, e-mail: pbeise@gmx.de, jmueller@uni-trier.de}}

\begin{document}
\maketitle

%

\begin{abstract}
It is known that, generically, Taylor series of functions holomorphic in the unit disc turn out to be universal series outside of the unit disc and in particular on the unit circle. Due to classical and recent results on the boundary behaviour of Taylor series, for functions in Hardy spaces and Bergman spaces the situation is essentially different. In this paper it is shown that in many respects these results are sharp in the sense that universality generically appears on maximal exceptional sets. As a main tool it is proved that the Taylor (backward) shift on certain Bergman spaces is mixing.   
\end{abstract}



\section{Introduction and main results}

For an open set $\om$ in $\C$ with $0 \in \om$ we denote by $H(\om)$ the (Fr{\'e}chet) space of functions holomorphic in $\om$ endowed with the topology of locally uniform convergence. Moreover, for  $f \in H(\om)$ we write $s_n f(z):=\sum_{\nu=0}^n a_\nu z^\nu$ for the $n$-th partial sum of the Taylor expansion $\sum_{\nu=0}^\infty a_\nu z^\nu$ of $f$ about $0$. A classical question in complex analysis is how the partial sums $s_n f$ behave outside the disc of convergence and in particular on the boundary of the disc. Based on Baire's category theorem, it can be shown that for functions $f$ in the unit disc $\DD$ generically the sequence $(s_n f)$ turns out to be universal outside of $\DD$. For precise definitions and a large number of corresponding results, we refer in particular to the expository article \cite{kahane}. For results on universal series in a more general framework see also \cite{BGENP}.
  
The situation changes if we consider classical Banach spaces of holomorphic functions. In a first part, we study the generic boundary behaviour of the  Taylor sections $s_n f$ of functions in Hardy spaces $H^p$ and in Bergman spaces $A^p$ of order $1\le p<\infty$. Subsequently we investigate the Taylor backward shift and the behaviour of Taylor sections in the case of Bergman spaces $A^p(\om)$ for more general domains $\om$.

Let $m$ denote the normalized arc length measure on the unit circle $\T$. For $1 \le p<\infty$, the Hardy space $(H^p, ||\cdot||_{p})$ is defined as the (Banach-) space of all $f \in H(\DD)$ such that 
\[
||f||_p:= \sup_{0 < r<1} \Big(\int_\T |f_r|^p dm  \Big)^{1/p}<\infty \;,
\]
where $f_r(z):=f(rz)$ for $z \in \DD$. For basic properties we refer e.g. to \cite{Du}. 

It is well known that each $f \in H^p$ has nontangential limits $f^*(z)$ at $m$-almost all $z$ on the unit circle $\T$ and that $f^* \in L^p(\T)$. Moreover, the mapping $f \mapsto f^*$ establishes an isometry between $H^p$ and the closure of the polynomials in $L^p(\T)$. As usual, we identify $f$ and $f^*$ and, in this way, $H^p$ and the corresponding closed subspace of $L^p(\T)$, which we for clarity also denote by $H^p(\T)$. In particular, the restrictions $s_n f|\T$ of the partial sums of the Taylor expansion of $f$ are the partial sums of the Fourier expansion of $f$.    
So it is consistent to write $s_n f$ also for the partial sums of the Fourier expansion $\sum_{\nu=-\infty}^\infty \hat{f}(\nu) z^\nu$ of $f \in L^1(\T)$.
  
According to the classical Carleson-Hunt theorem, for each $p>1$ and each $f \in L^p(\T)$ the partial sums $s_n f$ of the Fourier series converge to $f$ almost everywhere on $\T$.  Due to results of Kolmogorov, in the case $p=1$ we have convergence in measure and therefore, in particular, each subsequence of $(s_n f)$ has a subsubsequence 
converging almost everywhere to $f$. 

Our first result shows that, on the other hand, generically the partial sums turn out to have a "maximal" set of limit functions on closed sets of measure zero. 
We say that a property is satisfied for comeagre many elements of a complete metric space, if the property is satisfied on a residual set in the space. Moreover, for a compact subset $E$ of $\C$ we denote by $C(E)$ the space of continuous complex valued functions on $E$ endowed with the uniform norm $||\cdot||_{E}$.

\begin{thm} \label{Hardy}
Let $1 \le p < \infty$ and suppose $E$ to be a closed subset of $\T$ with vanishing arc length measure.  If $\Lambda \subset \N_0$ is infinite, then comeagre many $f$ in $H^p$ enjoy the property that for each $g \in C(E)$ a subsequence of $(s_n f)_{n \in \Lambda}$ tends to g uniformly on $E$. The same is true for comeagre many $f \in L^p(\T)$.   
\end{thm}  
The proof is given in Section \ref{sec2}.
\begin{remark}
Let $(a_\nu)_{\nu \in \Z}$ be an arbitrary sequence in $\C$ and let $s_n(z):=\sum_{\nu=-n}^n a_\nu z^\nu$ for $z \in \T$. As a consequence of Rogosinski summability (see \cite[p. 113]{Zyg}) we obtain the following result: If the sequence $(s_n(\zeta))$ is Ces{\`a}ro-summable  to $s$ at the point $\zeta \in \T$ and if a subsequence $(s_{n_j})_j$ of $(s_n)$ converges to some function $h$ uniformly on the closed set $\zeta E_1$, where $E_1=\{e^{\pm \pi i/(2k)}: k\in \N\} \cup\{1\}$, then necessarily $h(\zeta)=s$.

Since the Fourier series of a function $f \in L^1(\T)$ is Ces{\`a}ro-summable to $f(\zeta)$ at each point of continuity of $f$, this shows that the    
situation changes drastically if we consider the space $C(\T)$ of continuous functions on $\T$ instead of $L^p(\T)$ (or the disc algebra instead of $H^p$). While, according to the Kahne-Katznelson theorem (see e.g. \cite[p. 58]{Katz}, \cite{kahane}), each set of vanishing arc length is a set of divergence for $C(\T)$, the above considerations show that, for example, each uniform limit function $h$ of a subsequence of $(s_n f)$ on the closed set $\zeta E_1$ (having the single accumulation point $\zeta$) necessarily has to satisfy  $h(\zeta)=	f(\zeta)$. In particular, if $E \supset \zeta E_1$, for some $\zeta$, a maximal set of limit functions on $E$ as in Theorem \ref{Hardy} is not possible for functions in $C(\T)$.  On the other hand, pointwise universal divergence on arbitrary countable sets $E$ does hold generically (see \cite{HK} and \cite{Mue}). Further interesting results in this direction are found in the recent paper \cite{Papa^2}.
\end{remark}     

We focus now on the question of possible limit functions of $(s_n f)$  on parts of the unit circle for functions $f$ in the Bergman spaces $A^p$. 

Let $m_2$ denote the normalized area measure on $\DD$. For $1 \le p<\infty$, the Bergman space $(A^p, ||\cdot||_{p})$ is defined as the (Banach-) space of all $f \in H(\DD)$ such that 
\[
||f||_p:= \Big(\int_\DD |f|^p dm_2  \Big)^{1/p}<\infty. 
\]
For basic properties we refer to \cite{DS} and \cite{HKZ}. 
It is known (see \cite[p. 85]{DS}) that for $1 \le p <\infty$ and $f \in A^p$ the coefficients $a_n$ satisfy the condition
\[
a_n=o(n^{1/p}) \qquad (n \to \infty).
\]
The estimate combined with a result of Shkarin (see \cite{Shk}) implies that for $f\in A^p$ at most one continuous (pointwise) limit function can exist on each nontrivial subarc of $\T$. We shall show, in contrast, that maximal sets of limit functions generically exist on metrically large subsets of $\T$. A trigonometric  (or power) series on $\T$ is called universal in the sense of Menshov if each measurable function $g:\T \to \C$ is the almost everywhere limit of a  subsequence of the partial sums (see e.g. \cite{kahane}).

\begin{thm}\label{Bergman}
For all $1 \le p <\infty$ comeagre many $f$ in $A^p$ turn out to be universal in the sense of Menshov, i.e. for each measurable function $g:\T \to \C$ a subsequence of the partial sums $(s_n f)$ tends to $g$ almost everywhere on $\T$.  
\end{thm}

Again, the proof is given in Section \ref{sec2}.\\

We consider Bergman spaces on more general domains. For $\Omega$ a domain in the complex plane $\C$ and $1 \le p<\infty$ let $A^p(\Omega)$ be the Bergman space of all functions $f$ holomorphic in $\om$ and satisfying 
\[
||f||_p:=||f||_{\om,p}:=\Big(\int_\Omega |f|^p \, d\lambda_2\Big)^{1/p} < \infty, 
\]
where $\lambda_2$ denotes the 2-dimensional Lebesgue measure.  Again, $(A^p(\om),||\cdot||_{p})$ is a Banach space (see, e.g. \cite{DS}, also for further properties). In the case $\om=\DD$ we recover $A^p$, up to normalisation of the integral. 

If $\om$ is bounded with $0\in \om$ we define  $T=T_{\om,p}: A^p(\Omega) \to A^p(\Omega)$ by
\[
Tf(z):= \frac{f(z)-f(0)}{z} \quad (z \not=0), \qquad Tf(0):=f'(0).
\] 
 If $f(z)=\sum_{\nu=0}^\infty a_\nu z^\nu$ then 
\[
Tf(z)=\sum_{\nu=0}^\infty a_{\nu+1} z^\nu
\]
for $|z|$ sufficiently small. 
We call $T$ the Taylor (backward) shift on $A^p(\om)$. Backward shifts are studied intensively on the classical Hardy spaces $H^p$ and Bergman spaces $A^p$ (see e.g. \cite{CiRo}, \cite{Ross}). 

By induction it is easily seen that
\begin{equation}\label{iterates}
T^{n+1}f(z)=\frac{f-s_{n}f(z)}{z^{n+1}} \quad (z \not=0), \qquad T^{n+1}f(0)=a_{n+1},
\end{equation}
for $n \in \N_0$, where, as above, $s_n f$ denotes the $n$-th partial sum of the Taylor expansion of $f$ about $0$. This shows that the behaviour of the iterates $T^n f$ is closely related to the behaviour of the sequence of partial sum $s_n f$. 

Our aim is to study the dynamics of $T$ on $A^p(\om)$ and to deduce results concerning the boundary behaviour of the partial sums $s_n f$, for generic $f \in A^p(\om)$. For notions from topological dynamics and linear dynamics used in the sequel we refer to \cite{BayMath} and \cite{GEandPeris}. In particular, if $X$ is Banach space, a continuous linear operator $T:X \to X$ is called topologically transitive if for each pair of non-empty open sets $U,V$ in $X$ a positive integer $n$ exist with $T^n(U) \cap V \not= \emptyset$. If this condition holds for all sufficiently large $n$ (depending on $U,V$), then $T$ is said to be mixing. 
 It is known that $T$ is mixing on $A^2$ (\cite[p. 96]{GEandPeris}, see also \cite{GS}). Moreover, in \cite{BMM_X} it is shown that $T$ is mixing on $H(\om)$ for arbitrary open sets $\om$ with $0\in \om$ and having the property that each connected component of $\C_\infty \setminus \om$ (with $\C_\infty$ denoting the extended plane) meets $\T$. \\

We consider Carath{\'e}odory domains $\om$, that is, bounded simply connected domains whose boundary equals the outer boundary (see e.g. \cite[p. 171]{Con}). In addition, we suppose that $\overline{\om}$ does not separate the plane. If $\om$ is a Jordan domain, then both conditions are satisfied. 
  
\begin{thm}\label{hypShift}
Let $\Omega$ be a Carath{\'e}odory domain with $0 \in \Omega$ and so that $\overline{\om}$ does not separate the plane. If $\T \setminus \om$ contains some arc then $T$ is mixing on $A^p(\om)$ for all $1 \le p <\infty$. 
\end{thm}

For the proof we refer to Section \ref{sec3}.\\

Let $\om,p$ be so that $T$ is mixing on $A^p(\om)$. If $\Lambda \subset \N$ is infinite, then the Universality Criterion (\cite[Theorem 1]{GE} or \cite[Theorem 1.57]{GEandPeris}) shows that comeagre many $f \in A^p(\om)$ are universal for $(T^{n+1})_{n \in \Lambda}$, i.e. for comeagre many $f \in A^p(\om)$, we have that $\{T^{n+1}f: n \in \Lambda\}$ is dense in $A^p(\om)$.
From \eqref{iterates} we obtain 
\[
|T^{n+1}f|\ge |f-s_{n}f|
\]
on $\overline{\mathbb{D}}\cap \Omega$ for all $f \in H(\om)$.  
Theorem \ref{hypShift} immediately implies
 
\begin{corollary}\label{limitf}
Let $1\le p <\infty$ and let $\om$ be a Carath{\'e}odory domain with $0 \in \Omega$ and so that $\overline{\om}$ does not separate the plane. Moreover, suppose that $ \T\setminus \om$ contains some arc. If $\Lambda \subset \N_0$ is infinite, then for comeagre many $f$ in $A^p(\om)$ there is a subsequence of $(s_n f)_{n \in \Lambda}$ tending to 
$f$ in $A^p(\om \cap \DD)$ and locally uniformly on $\om \cap \overline{\DD}$. 
\end{corollary}   
 
Indeed: Let $f$ be a universal function for $(T^{n+1})_{n \in \Lambda}$. Then there is a sequence $(n_j)$ in $\Lambda$ with $T^{n_j+1}f \to 0$ in $A^p(\om)$ as $j \to \infty$ and therefore, in particular, $s_{n_j} f \to f$ in $A^p(\om\cap \DD)$. Moreover, since convergence in $A^p(\om)$ implies local uniform convergence (see e.g. \cite[p. 8]{DS}), we also have    
$s_{n_j}f \to f$ $(j \to \infty)$
locally uniformly on $\overline{\DD} \cap \Omega$. 

\begin{remark}
The first assertion (and thus also the assertion of Theorem \ref{hypShift}) does no longer hold for general (bounded) simply connected domains $\om$: If $\om_0$ is a domain with $\om_0 \supset \om$ and $\lambda_2(\om_0\setminus \om)=0$ then each sequence of polynomials which converges in $A^p(\om)$ also converges in $A^p(\om_0)$. Hence, if $s_{n_j} f \to f$ in $A^p(\om)$ then $f$ extends to a function holomorphic in $\om_0$. If we consider, for instance, $\om$ to be the unit disc minus a radial slit, then convergence of a subsequence of $(s_n f)$ in $A^p(\om)$ is only possible if $f$ extends to a holomorphic function in $\DD$.
\end{remark}

We consider the case $\DD \subset \om$, in which the corollary concerns the behaviour of $(s_n f)$ on $\T$. In particular, for comeagre many $f \in A^p(\om)$, a subsequence of $(s_n f)$ tends to $f$ locally uniformly on $\T \cap \om$. 
An interesting question is whether there are (finite) limit functions different from $f$ on parts of $\T \cap \om$. 

Due to a classical result of Fatou and M. Riesz (see e.g. \cite[Chapter 11]{Re}), for each function $f$ holomorphic in a domain $\om$ with $\DD \subset \om \not=\DD$ and each closed arc $\Gamma$ on $\om \cap \T$, the partial sums $s_n f$ converge uniformly to $f$ on $\Gamma$ if $a_n$ tends to 0. 
This holds in particular for functions in $H^1$ on each closed arc of holomorphy (if such an arc exists).
As already mentioned above, for $1 \le p <\infty$ and $f \in A^p$ the coefficients $a_n$ satisfy the condition
$a_n=o(n^{1/p})$. In general, the exponent $1/p$ is best possible. So, the result of Fatou and M. Riesz does not apply here. 
On the other hand, without posing any conditions on the growth of $(a_n)$, a recent result of Gardiner and Manolaki (see \cite{GM}) shows that arbitrary functions $f$ holomorphic in $\DD$ have the following remarkable property: \\

\textit {Let $(s_{n_k}f)$ be an arbitrary subsequence of $(s_n f)$ converging to a (finite) limit function $h$ pointwise on a subset $S$ of $\T$. If $f$ has nontangential limits $f^*(\zeta)$ for $\zeta \in S$, then $h=f^*$ almost everywhere (with respect to arc length measure) on $S$}.  \\

In particular, the theorem proves the special attraction of the "right" boundary function as a limit function in the case that $f$ extends continuously to some subarc of $\T$.

The following result implies that, on small subsets of $\om \cap \T$, even for functions that belong to $A^p(\om)$, where $\om$ is as in Theorem \ref{hypShift}, a maximal set of uniform limit functions generically exists. 
We recall that a closed subset of $\T$ is called a Dirichlet set if a subsequence of $(z^n)$ tends to $1$ uniformly on $E$.   

\begin{thm}\label{Dirichlet}
Let $1\le p <\infty$ and let $\om$ be a Carath{\'e}odory domain with $ \DD \subset \Omega$ and so that $\overline{\om}$ does not separate the plane. Moreover, suppose that $\T \setminus\om$ contains some arc. If $E \subset \T\cap \om$ is a Dirichlet set then comeagre many $f \in A^p(\om)$ enjoy the property that for each $h \in C(E)$ a subsequence of $(s_{n} f)$ tends to $h$ uniformly on $E$. 
\end{thm}

The proof is found in Section \ref{sec3}. 
  
\begin{remark} 
1. It is easily seen that each finite set in $\T$ is a Dirichlet set. Moreover, it is known that Dirichlet sets cannot have positive arc length measure (as also follows from the above results), but can have Hausdorff dimension $1$ (see e.g. \cite{kahane}). 

2. Let $f$ be holomrphic in $\DD$. It is known that the condition $a_n=o(n)$ implies that $(s_n f)$ is Ces{\`a}ro summable at each point $\zeta \in \T$ at which $f$ has an unrestricted limit (see e.g. \cite{Of}).  This holds in particular for functions in $A^p(\om)$ and all $\zeta \in \om \cap \T$ . Again, using results on Rogosinski summability, it can be shown (cf. \cite[Corollary 3.3]{KNP}) that in the case of the existence of a function $f \in A^p(\om)$ with $(s_{n_j}f)_j$ tending to $f+1$ uniformly on a compact set $E \subset \om \cap \T$, the set $E$ necessarily has to satisfy the following (Dirichlet type) condition at each point $z \in E$: For all sequences $(z_n)$ in $E$ with $z_n/z = 1+O(1/n)$ the sequence $(z_{n_j}/z)^{n_j}$ tends to 1 as $j \to \infty$. The condition is obviously satisfied for some subsequence $(n_j)$ of the positive integers if $E$ is a Dirichlet set. On the other hand, for $\zeta \in \T$ the set $\zeta E_N$, where $E_N:=\{e^{\pi i/k}: k \in \N,\, k \ge N\} \cup\{1\}$, does not satisfy the above condition for any $(n_j)$  at the (single) accumulation point $\zeta$. Thus, in particular, the assertion of Theorem \ref{Dirichlet} does not hold for compact sets $E\subset \om \cap \T$ containing some $\zeta E_N$.   
\end{remark}  

\section{Proofs of Theorems \ref{Hardy} and \ref{Bergman}}\label{sec2} 

The proofs are based on results on simultaneous approximation by (trigonometric or algebraic)  polynomials. For the case (algebraic) case of the Hardy space, this goes back to Havin (\cite{Ha}, see also \cite{Hr}). The general approach is inspired by and based on results of \cite{Hr}. 

 We consider a Banach space $X=(X, ||\cdot ||_X)$ with $X \subset L^1(\sigma)$ for some Borel set $M \subset \C$ and some Borel measure $\sigma$ supported on $M$.  We say that $X$ is trigonometric if $0 \not\in M$ and the trigonometric polynomials (i.e. the span of the monomials $P_n$, where $P_n(z)=z^n$ for $z \not=0$ and $n \in \Z$) form a dense subspace of $X$ with $\limsup_{n \to \infty } ||P_n||_X^{1/n} \le 1$ and $\limsup_{n \to \infty } ||P_{-n}||_X^{1/n} \le 1$. Correspondingly, we say that $X$ is analytic, if a similar condition holds "one-sided", that is, the (algebraic) polynomials are dense in $X$ with  $\limsup_{n \to \infty } ||P_n||_X^{1/n} \le 1$. In both cases $X$ is  separable since the corresponding polynomials with (Gaussian) rational coefficients also form a dense subset. The spaces $L^p(\T)$ are trigonometric and the spaces $H^p(\T)$ are analytic (with $\sigma=m$ the normalized arc length measure). If $E \subset \T$ then $C(E)$ is trigonometric and, according to Mergelians's theorem, also analytic if $E$ is a proper subset of $\T$.
Moreover, the Bergman spaces $A^p$ are analytic with $\sigma=m_2$ the normalized area measure on $\DD$ (see, e.g. \cite[p. 30]{DS}).  

Let $X^*$ denote the (norm) dual of $X$. If $X$ is analytic we define the Cauchy transform $K_X :X^* \to H(\DD)$ with respect to $X$ by 
\[
K_X\Phi(w)= \sum_{\nu=0}^\infty \Phi(P_\nu) w^\nu \quad (w \in \DD,\, \Phi \in X^*).
\] 
If $X$ is trigonometric then we define $K_X: X^* \to H(\C_\infty \setminus \T)$ (where $\C_\infty$ denotes the extended plane) as above, for $w \in \DD$, and in $\C\setminus \overline{\DD}$ by 
\[
K_X\Phi(w)= \sum_{\nu=1}^\infty \Phi(P_{-\nu}) w^{-\nu} \quad (|w|>1,\, \Phi \in X^*)
\]   
(note that always $K_X \Phi$ vanishes at $\infty$). 

Since the corresponding polynomials form a dense set in $X$, the Hahn-Banach theorem implies that $K_X$ is injective. We write $X^{c*}$ for the range $K_X(X^*)$ of $K_X$ (in $H(\DD)$ or $H(\C_\infty \setminus \T)$), the so-called Cauchy dual of $X$.  Moreover, we write $X_1 \oplus X_2$ for the direct sum of two Fr{\'e}chet spaces $X_1$ and $X_2$ (cf. \cite[p. 36]{GEandPeris}).

\begin{lemma}\label{perp}
Let $X$ and $Y$ be two trigonometric spaces. Then $X^{c*} \cap Y^{c*}=\{0\}$ if and only if the pairs of the form $(P, P)$, where $P$ ranges over the set of trigonometric polynomials, form a dense set in the sum $X \oplus Y$. The same holds for analytic spaces and algebraic polynomials.  
\end{lemma}

\textit{Proof}. Consider a functional $(\Phi, \Psi) \in (X \oplus Y)^*= X^* \oplus Y^*$. Then we have 
\[
0=(\Phi, -\Psi)(P_n, P_n)=\Phi(P_n) - \Psi(P_n)
\]
for all $n \in \Z$ if and only if $K_X\Phi=K_Y (\Psi)$. 

If $X^{c*}\cap Y^{c*}=\{0\}$ then $K_X \Phi=K_Y \Psi=0$. Since $K_X$ and $K_Y$ are injective, we obtain that $(\Phi, \Psi)=(0,0)$. Then the denseness of the span of the $(P_n,P_n)$ follows from the Hahn-Banach theorem.

If, conversely, the span of $(P_n,P_n)$ is dense in $X \oplus Y$ and if $\Phi$ and $\Psi$ are so that $K_X \Phi=K_Y \Psi$, then the Hahn-Banach theorem implies that $(\Phi, -\Psi)=(0,0)$ and thus also $K_X \Phi=K_Y\Psi=0$. 

The proof is similar for the analytic case, with $\Z$ replaced by $\N_0$. 
\hfill $\Box$ \\  

In order to apply the lemma we need more information about the Cauchy transforms involved. For $E \subset \T$ closed, the Riesz representation theorem says that $(C(E))^*$ is isometrically isomorphic to the space of complex Borel measures supported on $E$ and endowed with the total variation norm. If we identify $\Phi$ with the corresponding Borel measure $\mu$, then $\mu(P_n)=\int_E \zeta^n d\mu(\zeta)$, for all $n \in \Z$, and thus the Cauchy transform of $\mu$ is given by
\[
K_{C(E)}\mu(w)= \int_E \frac{d\mu(\zeta)}{1-\zeta w}\quad (w \in \C \setminus \T).
\]
Similarly, according to $\Phi(f)=\int_\T f \overline{h} dm$  ($f \in L^p(\T)$) we may identify $\Phi \in (L^p(\T))^*$ with a unique function $h \in L^q(\T)$, where $q$ is the conjugate exponent of $p$ (i.e. $pq=p+q$ for $p>1$ and $q=\infty$ for $p=1$). From this it is seen that
\[
K_{L^p(\T)}h(w)= \int_\T \frac{\overline{h}(\zeta)}{1-\zeta w}\, dm(\zeta) \quad (w \in \C \setminus \T).
\]  
\textit{Proof of Theorem \ref{Hardy}}. We start with the (simpler) trigonometric case $L^p(\T)$ 
and consider the family $(s_n)_{n \in \Lambda}$ (more precisely $f \mapsto s_n f|E$) as a family of (continuous) linear mappings from $L^p(\T)$ to $C(E)$. As mentioned above, $C(E)$ is separable.  The Universality Criterion (see e.g. \cite[Theorem 1]{GE} or \cite[Theorem 1.57]{GEandPeris}) implies that it is sufficient -- and necessary -- to show that for each pair $(f, g) \in  L^p(\T) \oplus C(E)$ and each $\varepsilon >0$ there exist a trigonometric polynomial $P$ and an integer $n \in \Lambda$ so that $||f-P||_{p}<\varepsilon$ and $||g -s_n P||_{E}<\varepsilon$.  Since $s_n P=P$ for all trigonometric polynomials $P$ and all sufficiently large $n$ (depending on the degree of $P$), it is enough to show that the pairs of the form $(P, P)$, where $P$ ranges over the set of trigonometric polynomials, form a dense set in $L^p(\T) \oplus C(E)$. Due to Lemma \ref{perp}, it suffices to show that 
\[
(L^p(\T))^{c*} \cap (C(E))^{c*}=\{0\}.
\]
To this aim, we consider $h \in L^q(\T)$ and $\mu$ a complex Borel measure on $\T$ supported on $E$ with $K_{L^p(\T)}h=K_{C(E)}\mu$. Then the measure $\nu:=\mu-\overline{h}m$ satisfies
\[
\nu(P_n)=\int_E \zeta^n d\mu(\zeta)-\int_\T \zeta^n \overline{h(\zeta)}\,dm(\zeta)=0 \quad (n \in \Z).
\]  
Since $C(\T)$ is trigonometric, we obtain $\nu=0$ and thus $\mu=\overline{h}m$. On the other hand, since $m(E)=0$, the measure $\mu$ is singular with respect to $m$. This shows that $\mu=0$ and then also $K_{L^p(\T)}h=0$ (actually $h=0$).  

The arguments are similar in the analytic case $H^p(\T)$ (cf. \cite{Hr}). Since $H^p(\T)$ is a closed subspace of $L^p(\T)$, the Hahn Banach theorem shows that again each functional $\Phi$ on $H^p(\T)$ is induced by some $h \in L^q(\T)$ (now, however, not in a unique way). We consider the Cauchy transforms $K_{H^p(\T)}\Phi$ and $K_{C(E)}\mu$ on $\DD$.  If $K_{H^p(\T)}\Phi=K_{C(E)}\mu$ then we obtain as above $\nu(P_n)=0$, now for all $n \in \N_0$. The F. and M. Riesz theorem then implies that $\nu$ is absolutely continuous with respect to $m$ and therefore the same is true for $\mu$. Still $\mu$ is also singular with respect to $m$. So again we obtain $\mu=0$ and then also $K_{H^p(\T)} \Phi=0$. \hfill $\Box$   \\

We now turn towards the proof of Theorem \ref{Dirichlet} As above, basically we need a result on simultaneous approximation. The corresponding deep considerations are found in \cite{Hr}. They complete former work of Kegejan and Talaljan (cf. \cite{Hr}). 

 Let $E$ be a proper closed subset of $\T$ with $m(E)>0$.  Then $E$ is said to satisfy Carleson's condition if 
\[
\ell(E):=\sum_{k} m(B_k) \log(1/m(B_k)) < \infty
\]
where $\T \setminus E=\bigcup_k B_k$ is the finite or countable union of the pairwise disjoint open arcs $B_k$. 
With this notation we have the following result. 

\begin{thm}\label{universal}
Let $1 \le p<\infty$ and $E \subset \T$ be closed with either $m(E)>0$ and $E$ not containing a closed subset of positive measure satisfying Carlesons's condition or else $m(E)=0$. If $\Lambda \subset \N_0$ is infinite, then comeagre many $f$ in $A^p$ enjoy the property that for each $g \in C(E)$ a subsequence of $(s_nf)_{n \in\Lambda}$ tends to $g$ uniformly on $E$.  
\end{thm}

\textit{Proof}. Let $1\le p < \infty$ be fixed. As in the proof of Theorem \ref{Hardy} it suffices to show that the pairs of the form $(P, P)$, where $P$ ranges over the set of polynomials, form a dense set in $A^p \oplus C(E)$. 

If $f \in H(\DD)$ and $0 \le r<1$ we write $M(r,f):=\max_{|z|\le r}|f(z)|$. With that, for $s>0$  we consider the (Banach) space $B_s$ of functions $f \in H(\DD)$ satisfying
\[
M(r,f)(1-r)^s \to 0 \qquad(r \to 1^-),
\]
equipped with the norm $||f||_{B_s}:=\max_{0 \le r <1}M(r,f)(1-r)^s$ (cf. \cite{Hr}).

The fundamental Theorem 4.1 in \cite{Hr} shows that for all $s>0$ there is a sequence of polynomials $Q_n$ with $Q_n \to 1$ in $B_s$ and $Q_n \to 0$ uniformly on $E$, as $n \to \infty$.  It is easily seen that $B_s$ is continuously embedded into $A^p$ for all $s <1/p$ (cf. \cite[pp. 78]{DS}). Thus, if we choose $s<1/p$, then we also have $Q_n \to 1$ in $A^p$. 

Let $D_q$ denote the Dirichlet space of order $q \ge 1$, that is, the space of all $f \in  H(\DD)$ with $f' \in A_q$. It is known that, for $p>1$, the Cauchy dual $(A^p)^{c *}$ of $A^p$ equals $D_q$, with $q$ the conjugate exponent (cf. \cite{HP}, \cite{CiRo}). Since the multiplication operator $f \mapsto P_1f$ is continuous on $A^p$, Theorem 1.3 of \cite{Hr} implies the assertion (actually for this we only need that the Cauchy dual contains all polynomials).  \hfill $\Box$  \\   

\textit{Proof of Theorem \ref{Bergman}}. It can be shown that for each $\varepsilon >0$ there is a closed set $E$ as in Theorem \ref{universal} and so that $m(\T \setminus E)< \varepsilon$. More explicitly, for given $0<\varepsilon<1$, let $N \in \N$ be so that $\sum_{j=0}^\infty (N+j)^{-2} < \varepsilon$. For such $N$ we consider $E_N=\bigcap_{j \in \N_0} E_{N,j}$ to be a Cantor set, where $E_{N,j}$ is defined by successive "cancellation" of $2^j$ open arcs of length $m_{N,j}:=2^{-j}(N+j)^{-2}$ (cf. \cite[p. 163]{Hr}). Then we obtain
\[
m(\T \setminus E_N) = \sum_{j=0}^\infty 2^{j} m_{N,j} = \sum_{j=0}^\infty \frac{1}{(N+j)^2} < \varepsilon
\]
and 
\begin{eqnarray*}
\ell(E_N) & = & \sum_{j=0}^\infty 2^{j} m_{N,j} \log(1/m_{N,j}) \\
& \ge &\sum_{j=0}^\infty 2^{j} m_{N,j} \log(2^j) =  \log(2) \sum_{j=0}^\infty \frac{j}{(N+j)^2} =\infty,  
\end{eqnarray*}
so that $E_N$ does not satisfy Carleson's condition. But then also no compact subset of positive measure satisfies the condition (see \cite[Theorem 5.1]{Hr}).

The proof of Theorem \ref{Bergman} now follows from Lusin's theorem by a diagonal argument where we choose an increasing sequence $E_N$ as above with $m(E_N) \to 1$ as $N \to \infty$ (cf. \cite{kahane}).

\begin{remark}
Let $B_0$ denote the little Bloch space, that is, the set functions $f \in H(\DD)$ satisfying
\[
M(r,f')(1-r) \to 0 \quad (r\to 1^-). 
\]
It is known that $B_0$ coincides with the closure of the polynomials in the Bloch space $B$ (and is thus in particular normed in that way), and that $B_0$ is contained in all $B_s$ for $s>0$ -- and therefore also in all $A^p$. For functions in $B_0$, the Taylor coefficients $a_n$ tend to $0$ (\cite[p. 80]{DS}). Moreover, one can show that the Cauchy dual of $B_0$ equals $D_1$ (see \cite{ACP}). From the considerations in the proof of Theorem \ref{universal} and Lemma \ref{perp} it is seen that $D_q \cap (C(E))^{c*}=\{0\}$ for all $q>1$ and all $E$ as in Theorem \ref{universal}. An interesting question is, which conditions on $E$ would guarantee $D_1 \cap (C(E))^{c*}=\{0\}$. The corresponding sets $E$ again turn out to be sets on which the partial sums $s_n f$, for generic $f \in B_0$,  have  maximal set $C(E)$ of uniform limit functions. In particular, it would be interesting to know if functions in $B_0$ can be universal in the sense of Menshov.    
\end{remark}
\section{Proofs of Theorems \ref{hypShift} and \ref{Dirichlet}}\label{sec3}


In order to see how tools from linear dynamics enter, we first give a short proof of Theorem \ref{hypShift} for the case of $\om$ being the unit disk. 

\begin{proposition}
Let $1 \le p<\infty$. Then $T$ is mixing on $A^p$.
\end{proposition}

The proof is a straight forward application of Kitai's criterion using the fact that $S:A^p \to A^p$, defined by $Sg(z)=zg(z)$, is a right inverse of $T$. Indeed: Lebesgue's dominated convergence theorem shows that $S^n g \to 0$ in $A^p$, for all $g \in A^p$. Moreover, \eqref{iterates} implies that $T^n p$ eventually vanishes for each polynomial $p$. Since the polynomials are dense in $A^p$ (see e.g. \cite[Theorem 3]{DS}), Kitai's criterion (see \cite[Theorem 3.4]{GEandPeris}) implies that $T:A^p \to A^p$ is mixing. \\

\begin{remark}
The operator $T$ is no longer mixing on the little Bloch space $B_0$ and on the Hardy spaces $H^p$, since in these cases the Taylor coefficients $a_n$ of all $f$ tend to $0$. But then it is easily seen that, for all $f$, the sequence $(T^n f)$ also tends to $0$ locally uniformly on $\DD$. This implies that $T$ cannot even be topologically transitive.        
\end{remark}

 For $M \subset \C_\infty$ we write
\[
M^\prime:=1/(\C_\infty\setminus \Omega)
\]
(with $1/\infty :=0$ and $1/0:=\infty$). Then for open sets $\om$ in $\C_\infty$ with $0 \in \om$  the set $\om'$
is a compact plane set. 

Let in the sequel $\om$ be a Carath{\'e}odory domain. It is readily seen that the Cauchy kernel provides a family of eigenfunctions for $T$. More precisely, for $\alpha \in \C$ we define 
\[
\gamma_{\alpha}(z)=1/(1-z \alpha )  \quad (z\in \C \setminus \{1/\alpha\}).
\]
Then, for each $\alpha$ in the interior of $\om'$, the function $\gamma_\alpha$ belongs to $A^p(\om)$ and $\gamma_\alpha$ is an eigenfunction for $T$ corresponding to the eigenvalue $\alpha$. In particular, since $\om'$ coincides with the closure of its interior, the compact set $\om'$ is contained in the spectrum of $T$. On the other hand, one observes that in case $\alpha \in 1/\om$   
\[
S_\alpha g(z):=\frac{zg(z)- g(1/\alpha)/\alpha}{1-z\alpha} 
\]
(continuously extended at the point $1/\alpha$) defines the continuous inverse operator to $T-\alpha I$ (with $I$ being the identity operator on $A^p(\om)$). 
This shows that the spectrum of $T$ on $A^p(\om)$ equals $\om'$. 
In the case $p<2$ the functions $\gamma_\alpha$ belong to $A^p$ also for $\alpha \in \partial (1/\om)=\partial (\om')$, and therefore they are also eigenfunctions.  In particular, the spectrum equals the point spectrum in that case. 

It is known that a sufficient supply of unimodular eigenvalues implies that $T$ is topologically transitive or even mixing (see e.g. \cite{BayMath}, \cite{GEandPeris}).

We are interested mainly in the case $\DD \subset \om$. Then unimodular eigenvalues exist only for $p<2$.
Therefore, in the case $p\ge 2$ an approach to universality properties via unimodular eigenvalues is no longer possible. Instead, for $p\ge 2$
we consider certain "integral means" of  eigenvectors  corresponding to unimodular eigenvalues for $p<2$:

Let $\Gamma\subset \T$ be a closed arc.
We consider the Cauchy integral $f \in H(\C \setminus \Gamma)$, defined by
\[
f(w) =\int_\Gamma \frac{d\zeta}{\zeta-w}   \qquad (w \not\in \Gamma).
\]
It is well known (see e.g. \cite[Theorem 1.7]{HKZ}) that 
\[
\int_\T \frac{dm(\zeta)}{|\zeta-w|}= O \left(\log\frac{1}{1-|w|}\right) \qquad(|w| \to 1^-),
\]
which implies that $|f|^p$ is locally integrable on $\C$ for all $1 \le p <\infty$ and thus, in particular, $f \in A^p(\om)$. 

For $\alpha \in \C$ we define $f_\alpha=f_{\alpha,\Gamma} \in H(\C \setminus \alpha^{-1}\Gamma)$ (with $\infty \Gamma :=\emptyset$) by 
\[
f_\alpha (z):=f(\alpha z)=\int_\Gamma \frac{d\zeta}{\zeta-\alpha z} \qquad(z \not\in \alpha^{-1}\Gamma).
\]
We consider $\Gamma, A \subset \T$ to be  closed arcs with $A^{-1}\Gamma \subset \T \setminus \overline{\om}$. 
 From the above considerations it follows that $f_\alpha \in A^p(\om)$ for all $\alpha \in A$ and all $1\le p < \infty$.

\begin{lemma}\label{dense sub}
Let $\om$ be a Carath{\'e}odory domain so that  $\overline{\om}$ does not separate the plane. If $\T \setminus \om$ contains an arc, then closed arcs $A\subset \T$ and $\Gamma\subset \T$ exist with the property that for each subset $B$ of $A$ such that the closure in $A$ has positive $m$-measure the span of $\{f_\alpha =f_{\alpha,\Gamma}: \alpha \in B\}$ is dense in $A^p(\om)$, for all $1\le p <\infty$.
\end{lemma}

\textit{Proof}.
According to the Farrell-Markushevich theorem, the set of polynomials, and thus in particular $A^p(\om)$ for $p>1$, is dense in $A^1(\om)$ (see e.g. \cite[p. 173]{Con}). Therefore, it suffices to prove the result for $p>1$. 

Since $\T\setminus \om$ contains an arc, there are closed arcs $A,\Gamma \subset \T$ such that $\text{dist}(A^{-1}\Gamma, \T \cap \om)>0$. Moreover, since $\overline{\om}$ does not separate the plane, the arc $\Gamma$ can be chosen so small that, in addition, the open set $\{\alpha \in \C :\text{dist}(\Gamma, \alpha \om)>0\}$ contains a connected open set $U$ with $0 \in U$ and so that $A \subset \partial U$. 
From the above considerations it is seen that for $\alpha \in U$ the function $f_\alpha$ is holomorphic in a neighbourhood of $\overline{\om}$ and thus $f_\alpha \in A^p(\om)$. 

Let $\Phi\in A^p(\om)^*$ be given. If $\alpha \in U$ is fixed then there are neighbourhoods $V$ of $\overline{\om}$ and $W$ of $\alpha$ so that $f_\beta$ is holomorphic in $V$ for all $\beta \in W$. This implies that 
\[
\frac{f_\beta(z) -f_\alpha(z)}{\beta-\alpha} \to z \int_\Gamma \frac{d\zeta}{(\zeta-\alpha z)^2} \qquad (\beta \to  \alpha) 
\]
uniformly on $\overline{\om}$ and therefore also in $A^p(\om)$. From this it is seen that  the function $h:U \cup A \to \C$, defined by $h(\alpha):=\Phi(f_\alpha)$, is holomorphic on $U$. It suffices to show that $h$ vanishes identically if it vanishes on the set $B$.
(Indeed: In this case $h^{(\nu)}(0)=\Phi(P_\nu)\, \int_\Gamma 1/\zeta^\nu\, d\zeta=0$ for all $\nu\geq 0$ so that, again by the Farrell-Markushevich theorem, $\Phi=0$. The assertion then follows by the Hahn-Banach theorem.) 

From the definition of $f_\alpha$ it is seen that there exists a neighbourhood $D$ of $A$ relative to $\overline{\DD}$ such that $\{f_\alpha: \alpha \in D\}$ is a bounded family in $A^p(\om)$. We can choose $D$ in such a way that the interior $D^o$ (with respect to $\C$) is a Jordan domain with piecewise smooth boundary (a sector, for instance).  

If $\alpha \in A$ and if $(\alpha_n)$ is a sequence in $D$ with $\alpha_n \to \alpha$, then 
$f_{\alpha_n} \to f_\alpha$ pointwise on $\om$.  
Since $p>1$, the boundedness of the family $\{f_\alpha: \alpha \in D\}$ implies that $h(\alpha_n) \to h(\alpha)$ (see Lemma 1.10 of \cite{Con}). This shows that $h$ is continuous on $A$. Thus, if $h|_B=0$ we have vanishing nontangential limits of $h$ at all points of the closure of $B$ in $A$. Moreover, the boundedness of the family $\{f_\alpha: \alpha \in D\}$ implies the boundedness of $h$ on $D$. Since the closure of $B$ in $A$ has positive measure, we obtain that $h=0$ on $D^o$ (cf. \cite[Theorem 10.3]{Du}) and then also on $U$. 
\hfill $\Box$\\

\begin{remark} \label{oldcase}
The proof of the Lemma shows that in the case $p<2$, for each arc $A \subset \om'$ and each subset $B$ of $A$ such that the closure in $A$ has positive measure, the span of $\{\gamma_\alpha: \alpha \in B\}$ is dense in $A^p(\om)$. The corresponding approximation result appears for $p=1$ as a special case of the main theorem from \cite{Be}. 
 \end{remark}

\textit{Proof of Theorem \ref{hypShift}}.
Let $A$ and $\Gamma$ be closed arcs of $\T$ as in Lemma \ref{dense sub}.Then the span of $\{f_\alpha : \alpha \in A\}$ is dense in $A^p(\om)$. For $f\in H(\om)$ we have 
\[
T^nf(z)= \frac{1}{2\pi i} \int_\gamma \frac{f(\xi)}{\xi^n\, (\xi-z)}\, d\xi  
\]
with $\gamma$ some closed path with $\ind_\gamma(0)=\ind_\gamma (z)=1 $ and $\ind_\gamma(w) =0 $ for all $w\notin \om$. Applying the Cauchy formula to the functions $f_\alpha$, this yields
\[
T^n f_\alpha(z)= \alpha^n \int_\Gamma \frac{d\zeta}{(\zeta-\alpha z)\, \zeta^n},
\]
for $\alpha \in A$. By partial integration we get
\[
 \int_\Gamma \frac{d\zeta}{(\zeta-\alpha z)\, \zeta^n} =\frac{1}{n-1}\, \left( \int_\Gamma \frac{-d\zeta}{(\zeta-\alpha z)^2\, \zeta^{n-1}} - \frac{1}{(b-\alpha z)\, b^{n-1}}+\frac{1}{(a-\alpha z)\, a^{n-1}}\right)
\]
with $a$ and $b$ the endpoints of $\Gamma$. Hence, for every open neighborhood $U$ of $\alpha^{-1}\Gamma$, we have that $T^n f_\alpha(z) \rightarrow 0$ as $n$ tends to $\infty$ uniformly on $\om\setminus U$. Since 
\[
|T^{n}f_\alpha (z)| \le 2\pi \int_\T \frac{dm(\zeta)}{|\zeta-\alpha z|} \qquad(z \in \DD)
\]
and since $U$ can be chosen of arbitrary small area, this shows that $||T^n f_\alpha||_p\rightarrow 0$ for $n\to \infty$. If we define $S_n$ on the span of $\{f_\alpha : \alpha \in A\}$ by 
\[
S_n f_\alpha(z):=\frac{1}{\alpha^n} \int_\Gamma \frac{\zeta^n}{\zeta-\alpha z} \, d\zeta
\] 
(and linearity) we have $T^nS_n f_\alpha =f_\alpha$. Moreover, $||S_n f_\alpha||_p\rightarrow 0$ for $n\rightarrow \infty$ follows by similar arguments as above. An application of Kitai's criterion (cf. \cite[Theorem 3.4]{GEandPeris}) yields the assertion. 
\hfill $\Box$\\

\begin{remark} \label{comeagre}
In the case $p<2$, Theorem \ref{hypShift} may also be deduced from Remark \ref{oldcase} and \cite[Theorem 5.41]{BayMath}.
Moreover, in this case, for $\om$ is as in Theorem, the operator $T$ is also chaotic on $A^p(\om)$. Indeed: Let $A$ be an arc in $\T \cap \om'$. Since the span of $\{\gamma_\alpha: \alpha \in A, \alpha \textrm{ a root of unity}\}$ consists of periodic points, Lemma \ref{dense sub} (cf. Remark \ref{oldcase}) implies that the periodic points are dense in $A^p(\om)$. 
This is no longer true for $p \ge 2$ and $\DD \subset \om$, in which case actually no periodic points exist (cf. \cite[p. 96]{GEandPeris}).  
\end{remark}	


\textit{Proof of Theorem \ref{Dirichlet}}. Let $\Lambda \subset \N_0$ be infinite with $z^{n+1} \to 1$ uniformly on $E$ as $n \to \infty$, $n \in \Lambda$. 
According to Mergelian's theorem (note that $E$ has connected complement), the polynomials are dense in $C(E)$. So we can assume $h \in A^p(\om)$. 

Let $f$ be universal for $(T^{n+1})_{n \in \Lambda}$ (which is the case for comeagre many $f$ in $A^p(\om)$). Since convergence in $A^p(\om)$ implies local uniform convergence, there are $n_j$ in $\Lambda$ with 
$T^{n_{j}+1} f \to f-h$ $(j \to \infty)$ locally uniformly on $\om$ and thus in particular uniformly on $E$. Then also 
\[
z^{n_{j}+1} T^{n_{j}+1}f(z) \to (f-h)(z) \qquad (j \to \infty)
\]
uniformly on $E$ and therefore
\[
s_{n_{j}}f(z) =f(z)-z^{n_{j}+1} T^{n_{j}+1}f(z) \to h(z)  \qquad (j \to \infty)
\]
uniformly on $E$.\hfill $\Box$\\

\end{document}